\documentclass[a4paper,twoside,11pt]{article}
\usepackage{amssymb}
\usepackage{amsmath}
\usepackage{amsthm}
\usepackage{latexsym}
\usepackage{amsfonts}
\usepackage{graphicx}
\usepackage{amscd}
\usepackage{graphics}

\newcommand{\dis}{\displaystyle}
\theoremstyle{plain}
\newtheorem{thm}{Theorem}[]   
\newtheorem{prop}[thm]{Proposition}
\newtheorem{lem}[thm]{Lemma}

\theoremstyle{definition}

\newtheorem*{Proof}{Proof}

\newcommand{\bbb}[1]{\mbox{\boldmath$#1$}}

\newcommand{\el}{\ell}

\newcommand{\ra}{\;\rightarrow\;}

\newcommand{\al}{\alpha}
\newcommand{\bi}{\beta}
\newcommand{\ga}{\gamma }

\newcommand{\de}{\delta }

\newcommand{\e}{\varepsilon }

\newcommand{\thi}{\theta }
\newcommand{\vthi}{\vartheta }

\newcommand{\la}{\lambda }

\newcommand{\si}{\sigma }

\newcommand{\oo}{\omega}
\newcommand{\C}{\mathbb{C}}

\newcommand{\R}{\mathbb{R}}

\newcommand{\ssum}{\sum\limits}

\newcommand{\ld}{\ldots}

\newcommand{\sm}{\smallsetminus}

\newcommand{\qb}{$\quad\blacksquare$}

\begin{document}
\pagestyle{myheadings}
\markboth{V. Nestoridis`}{V. Nestoridis}
\title{\bf A project about chains of spaces,\\ regarding topological and algebraic genericity and spaceability}
%
%
\author{V. Nestoridis\\
National and Kapodistrian University of Athens}
\date{}
\maketitle
\begin{abstract}
We present the example of $\el^p$ spaces, where we examine results of topological and algebraic genericity and spaceability.
 At the end of the paper we include a project with other chains of spaces, mainly of holomorphic functions, where a similar
  investigation can be done.
\end{abstract}
{\em AMS classification numbers}: primary, 15A03, 46E10, 46E15, secondary, 30H10, 30H20, 30H35.\smallskip
\setcounter{section}{-1}
\section{Introduction}\label{sec1}
\noindent
We consider chains of spaces $X_i\subset X_j$, $X_i\neq X_j$ for $i<y$ and we examine if $X_i$ is an $F_\si$
meager set in $X_j$ (topological genericity); equivalently, if $X_j-X_i$ is a $G_\de$-dense subset of $X_j$.
The main tool towards this is Baire's Category theorem for complete metric spaces. Furthermore,
we examine if $X_j-X_i$ contains a vector space, except 0, dense in $X_j$ (algebraic genericity). Finally,
we examine if $X_j-X_i$ contains a vector space, except 0, which is infinite dimensional and closed in $X_j$ (spaceability).
 One can also examine, if $\dis\bigcup_{i<j}X_i$ is an $F_\si$ meager in $X_j$ (topological genericity) and
 if $X_j\sm\Big(\dis\bigcup_{i<j}X_i\Big)$ contains vector spaces except 0 (algebraic genericity and spaceability).

In Sections \ref{sec1}, \ref{sec2} and \ref{sec3} we treat the example of $\el^p$ spaces. We do not include
the study of $\el^p\sm\dis\bigcup_{q<p}\el^q$ which follows easily in a similar way as the cases included in
Sections \ref{sec1},\ref{sec2},\ref{sec3}. Most, but not all, of the examples, concerning $\el^p$ spaces,
treated in Sections \ref{sec1},\ref{sec2},\ref{sec3} can be found in \cite{1} and \cite{2}, where
the algebraic genericity and spaceability for $\el^p\sm\dis\bigcup_{q<\bi}\el^p$ $\bi\le p$ are also included.

For algebraic genericity and spaceability we are also  referring to the works mentioned in \cite{1}, \cite{2}
and especially these
of Aron, Gurariy and Seoane\,-\,Sep\'{u}velda, Bernal\,-\,Gonzalez, F. Bayart and S. Charpentier.
For the topological genericity
we refer to \cite{3} below and the references there in, especially the works of J.-P. Kahare,
K.-G. Grosse-Erdmann and \cite{4}.

The chain of $\el^p$ spaces can be extended adding intersections of such spaces, as well as,
the inclusions $\el^p\subset c_0\subset\el^\infty\subset H(D)\subset\C^{N_0}$, where every
sequence $(a_n)$ in $\el^\infty$ can be identified with the function $f(z)=\sum a_nz^n$, which
is holomorphic in the open unit disc $D$ of the complex plane $C$. The cartesian product $\C^{N_0}$
can be identified with the set of formal power series $\sum a_nz^n$, $(a_n)\in C^{N_0}$.

In Section \ref{sec4} we present a project with other chains of spaces, where similar questions may
be investigated. They are mostly spaces of holomorphic functions in the disc $D$, but they also include
some sequence spaces, via the identification $(a_n)=\sum a_nz^n$. We can also add to the project other
spaces of holomorphic functions as Dirichlet spaces or Bloch spaces, as well as spaces of holomorphic
functions on domains in $\C^d$, $d\ge1$ and not only on $D$. In all previous cases the spaces are
complete metrizable topological vector spaces, in fact $F$-spaces, and the injections $X_i\subset X_j$
are continuous. What about spaces of harmonic functions on domains of $\R^{n2}$?

In Sections \ref{sec5} and \ref{sec6} we give two theorems concerning Hardy spaces on $D$:
Theorem \ref{thm8} about topological genericity and Theorem \ref{thm9} about algebraic genericity.
These results are in the frame of the project presented in Section \ref{sec4}. A result complementing
Theorem \ref{thm8} is the following (see also \cite{5} and \cite{6}). \vspace*{0.5cm} \\
\noindent
{\bf Theorem A}. {\em Let $0<p<+\infty$. Then there exists a holomorphic function $f$ on the open
unit disc $D$, such that 1) and 2) below hold.
\begin{enumerate}
\item[1)] $\sum_{0<r<1}\dis\int^{2\pi}_0|f(re^{i\thi})|^qd\thi<+\infty$ for all $0<q<p$
\item[2)] $\dis\sup_{0<r<1}\dis\int^b_a|f(re^{i\thi})|^\de md\thi=\infty$ for all $p\le\de<+\infty$ and $a<b$.
\end{enumerate}

The set of such functions $f$ is a $G_\de$ and dense subset of $\dis\bigcap_{q<p}H^q$ endowed with its natural topology}.

For the proof we use the fact that for $\oo\in\R$ and $\ga=\dfrac{1}{p}$ the function
$g(z)=\dfrac{1}{(z-e^{i\oo})^\ga}$ belongs to $H^q$ for all $q\in(0,p)$ but not to $H^p$.
The proof of the previous theorem is similar to that of Theorem \ref{thm8} and is omitted.

A result complementing Theorem \ref{thm9} is the following theorem, whose proof is similar
to that of Theorem \ref{thm9} and is omitted.\vspace*{0.5cm}\\
\noindent
{\bf Theorem B}. {\em Let $0<p\le q\le+\infty$ and $\al<\bi$ be fixed.
Then $\Big(\dis\bigcap_{\bi<p}H^\bi\sm H^q_{[a,b]}(D)\Big)\cup\{0\}$ contains a vector space
dense in $\dis\bigcap_{\bi<p}H^\bi$ endowed with it natural topology.

The space $\dis\bigcap_{\bi<p}H^\bi$ can be replaced by $H^p$ provided $p<q$}.
\section{Topological genericity for the $\el^p$ spaces}\label{sec1}
\noindent

We will deal with the following spaces
\[
\el^\infty\supset c_0\supset\bigcap_{p>b}\el^p\supset\el^b\supset\bigcap_{p>a}\el^p\supset
\el^a\supset\bigcap_{p>0}\el^p
\]
where $0<a<b<+\infty$.

All inclusions are strict. For instance \\
$(1,1,\ld)\in\el^\infty-c_0,\;\Big(\dfrac{1}{n^{\frac{1}{b+1}}}\Big)^\infty_{n=1}
\in c_0-\dis\bigcap_{p>b}\el^p$ ,\; $\Big(\dfrac{1}{n^{1/b}}\Big)^\infty_{n=1}\in\Big(\dis\bigcap_{p>b}\el^p \Big)\sm\el^b,\;
\Big(\dfrac{1}{n^{\frac{1}{\ga}}}\Big)^\infty_{n=1}\in$\\
$\el^b-\dis\bigcap_{p<a}\el^p$ with $\ga=\dfrac{a+b}{2}$,\; $\Big(\dfrac{1}{n^\frac{1}{a}}\Big)^\infty_{n=1}\in$
$\dis\bigcap_{p>a}\el^p\sm\el^a,\Big(\dfrac{1}{n^{\frac{1}{x}}}\Big)^\infty_{n=1}
\in\el^a\sm\dis\bigcap_{p>0}\el^p$ with $x=\dfrac{a}{2}$.

Also all above spaces are metrizable complete topological vector spaces, when endowed with their natural topologies.
The spaces $\el^\infty$, $c_0$, $\el^p$ with $1\le p<+\infty$ are Banach spaces. The space $\dis\bigcap_{p>b}\el^p$
with $1\le b<+\infty$ is a Fr\'{e}chet space. We consider a strictly decreasing sequence $p_m\downarrow b$ and the
distance in this space is defined by
\[
d(f,g)=\sum^\infty_{m=1}\frac{1}{2^m}\frac{\|f-g\|_{p_m}}{1+\|f-y\|_{p_m}}, \ \ \text{where} \ \ \|F\|_p=
\bigg(\sum^\infty_{n=1}|F(n)|^p\bigg)^{1/p}.
\]

For $0<p<1$ the space $\el^p$ is not a Banach space. It is a metrizable complete topological vector space
with metric $d_p(f,g)=\ssum^\infty_{n=1}|f(n)-g(n)|^p$. For $0\le a<1$ the space $\dis\bigcap_{p>a}\el^p$
is a metrizable complete topological space. Let $p_m\downarrow a$. Then,the metric in this space is
\[
d(f,g)=\sum^\infty_{m=1}\frac{1}{2^m}\frac{d_{p_m}(f,g)}{1+d_{p_m}(f,g)}
\ \ \text{where} \ \ d_p(f,g)=\sum^\infty_{n=1}|f(n)-g(n)|^p.
\]
Let $Y$ and $X$ two of the previously mentioned spaces with $X\subset Y$ and $Y-X\neq\emptyset$. Then,
the injection map $I:X\ra Y$, $I(f)=f$ is linear continuous and it is not sufjective. Then, according to a theorem
of Banach its image $I(X)=X$ is meager in $Y$; that is, $X$ is contained in a denumerable union of closed in $Y$
subsets with empty interior (in $Y$). In some cases we will show that $X$ is equal to such a set; that is, $X$
is an $F_\si$ meager subset of $Y$; equivalently $Y\sm X$ is a $G_\de$ and dense subset in $Y$, while in
the general case $Y-X$ is residual in $Y$.

Let $Y,X$ be two spaces among the previously mentioned ones, with $X\subset Y$ and $Y\sm X\neq\emptyset$.
If $X=c_0$, then $Y=\el^\infty$. In this case $X$ is a closed vector subspace of $Y$ and hence it has empty
interior in $Y$. Thus, $X=c_0$ is an $F_\si$ meager subset of $Y=\el^\infty$.

Next we consider the case $X=\el^p$ with $0<p<+\infty$.
\begin{prop}\label{prop1}
Let $X=\el^p$ with $0<p<+\infty$ and $Y\supset X$, $Y-X\neq\emptyset$ be one of the spaces mentioned
previously. Then $X$ is an $F_\si$ meager subset in $Y$.
\end{prop}
\begin{Proof}
We have that convergence in $Y$, $f_m\ra f$, as $m\ra+\infty$, implies pointwise convergence,
$f_m(n)\ra f(n)$ as $m\ra+\infty$ for all $n=1,2,\ld$. Since $Y-X\neq\emptyset$, there exists $g\in Y\sm X$.
Fix such a $g$. We have $X=\dis\bigcup^\infty_{M=1}E_M$, where $E_M=\Big\{f\in Y:\ssum^N_{n=1}|f(n)|^p\le M$
 for all $N=1,2,\ld\Big\}$. We will show that $E_M$ is closed in $Y$.

Indeed let $f^m\in E_M$ for $m=1,2,\ld$ and $f^m\xrightarrow[m\ra+\infty]{}f$ in $Y$. Since convergence in $Y$
implies pointwise convergence we have $\ssum^N_{n=1}|f^m(n)|^p\xrightarrow[m\ra+\infty]{}\ssum^N_{n=1}|f(n)|^p$
for every $N=1,2,\ld$. As $f^m\in E_M$ we have $\ssum^N_{n=1}|f^m(n)|^p\le M$. This implies
 $\ssum^N_{n=1}|f(n)|^p\le M$ for all $N=1,2,\ld$. Thus, $f\in E_M$ and we have proved that $E_M$ is
 closed in $Y$. Thus, $X$ is an $F_\si$ subset of $Y$. Since, according to the theorem of Banach, $X$ is meager
  in $Y$, the proof is completed. However, in order to avoid the use of the theorem of Banach, we will show
  that $E^\circ_M=\emptyset$. If not, there exists $f\in Y$ with $f\in E^\circ_M$. Let $g$ be in $Y-X$.
  Then, $f+\dfrac{1}{k}g\underset{k}{\ra}f$ in $Y$ because $Y$ is a topological vector space. Since $f$
  belongs to the open in $Y$  set $E^\circ_M$, there exists $k$ so that $f,f+\dfrac{1}{r}g\in E^\circ_M\subset E_M\subset X$.
  Since $X$ is also a vector space, it follows that $g\in X$ which contradicts the fact that $y\in Y\sm X$.

The proof is completed. \qb
\end{Proof}

The remaining case is $X=\dis\bigcap_{p>a}\el^p$, $0\le a<+\infty$. Then $Y\supset X$, $Y-X\neq\emptyset$ is
one of the previously mentioned spaces. Therefore, there exists $b\in(a,+\infty)$ such that $X\subset\el^b\subset Y$
and $Y-\el^b\neq\emptyset$. It follows that $\el^b$ is an $F_\si$ meager subset of $Y$, according to
Proposition \ref{prop1}. Thus, $X$ is contained in $\el^b$ which is a denumerable union of closed in $Y$ sets
with empty interiors. It follows that $X$ is meager in $Y$ and we arrived to this conclusion without using Banach's
theorem. Since $X=\dis\bigcap_m\el^{p_m}$ where $p_m\in(a,+\infty)$ is a strictly decreasing sequence
 $p_m\downarrow a$ and each $\el^{p_m}$ is an $F_\si$ meager subset of $Y$ we can not conclude
 that $X=\dis\bigcup_{p>a}\el^p$ is an $F_\si$ subset of $Y$. I do not know the answer if $\dis\bigcap_{p>a}\el^p$
 is an $F_\si$ in $Y$ or not.
\section{Algebraic genericity for the $\el^p$ spaces}\label{sec2}
\noindent

Let $X$ and $Y$ be two spaces mentioned above, such that $Y\supset X$ and $X\neq Y$. Thus,
 we say that there is algebraic genericity for the couple $(Y,X)$ if there is a vector subspace $F$ of $Y$
 dense in $Y$, such that $F\sm\{0\}\subset Y\sm X$.

If $Y=\el^\infty,1$ do not know the answer whether there is algebraic genericity for the couple $(\el^\infty,X)$
or not. May be the difficulty is that $\el^\infty$ is not separable. In all other cases we have algebraic genericity.
The essential lemma is the following.
\begin{lem}\label{lem2}
Let $0<b<+\infty$, $Y=\dis\bigcap_{p>b}\el^p$ and $X=\el^b$. Then we have algebraic genericity for the couple $(Y,X)$.
\end{lem}
\begin{Proof}
Let $x_j$, $j=1,2,\ld$ be an enumeration of all elements $x_j=f\in C_0$ with $f(n)\in Q+iQ$ for all $n=1,2,\ld$
and such that there exists $n_j\in N$ such that $x_j(n)=f(n)=0$ for all $n>n_j$. Then the set $\{x_j:j=1,2,\ld\}$
is dense in $Y$.

Let also $A_j$, $j=1,2,\ld$, be a sequence of infinite subsets of $\{1,2,\ld,\}$ pairwise disjoint. Since $A_j$ is
infinite denumerable, there is a function $y_j:A_j\ra\C$, such that,
\[
\sum_{n\in A_j}|y_j(n)|^b=+\infty \ \ \text{and} \ \ \sum_{n\in A_j}|y_j(n)|^p<+\infty \ \ \text{for all} \ \ p>b.
\]
We extend $y_j$ on $N$ by setting $y_j(n)=0$ for all $n\in N\sm A_j$. Thus, $y_j\in Y\sm X$. We consider the sets
\[
V_j=\bigg\{g\in Y:d_{p_k}(g,0)<\frac{1}{j} \ \ \text{for} \ \ k=1,2,\ld,j\bigg\},
\]
where $p_k$, $k=1,2,\ld$ is a strictly decreasing sequence converying to $b$ and
\[
d_p(g,0)=\bigg(\sum_{n\in N}|g(n)|^p\bigg)^{1/p} \ \ \text{for} \ \ p\ge1 \ \ \text{and} \ \ d_p(g,0)=
\sum_{n\in N}|g(n)|^p \ \ \text{if} \ \ 0<p<1.
\]
Then, the sequence $V_j$, $j=1,2,\ld$ is a base of neiborhoods of 0 in $Y=\dis\bigcap_{p>b}\el^p$. Since $Y$ is a
topological vector space we have $\dis\lim_{c\ra0}cy_j=0$ in $Y$. Thus, there exists $c_j\neq0$ with $c_jy_j\in V_j$.
We will construct a sequence $f_j\in Y=\dis\bigcap_{p>b}\el^p$, such that $f_j-x_j\in V_j$ for all $j=1,2,\ld$ and
any non-zero finite linear combination of the $f_j$  will belong to $Y-X$. Then, the linear spaces  $F$ of the $f_j$'s,
$j=1,2,\ld$ will be a vector subspace of $Y$ dense in $Y$, such that $F\sm\{0\}\subset Y-X$. That is, we have
algebraic genericity.

It suffices to set $f_j(n)=x_j(n)$ for $n\le n_j$, $f_j(n)=c_jy_j(n)$ for $n\in A_j-\{1,2,\ld,n_j\}$ and $f_j(n)=0$
otherwise (or simply $f_j=x_j+c_jy_j)$. One early checks that $f_j\in Y\sm X$ and that $f_j-x_j\in V_j$.
 Let $f=a_1f_1+\cdots+a_Nf_N$ with $a_N\neq0$. Since $Y$ is a vector space, it follows that $f\in Y$, Also
\[
\sum_{n\in A_N\atop n>n_N}|f(n)|^b=|a_N|^b|c_N|^b\sum_{n\in A_N\atop n>n_N}|y_N(n)|^b=+\infty.
\]
Therefore, $f\notin \el^b$. It follows that $f\in Y\sm X$.

The proof is completed.  \qb
\end{Proof}
\begin{prop}\label{prop3}
Let $Y$ and $X$ be two of the previously mentioned spaces, $Y\neq\el^\infty$, $Y\supset X$, $Y\neq X$.
Then there is algebraic genericity for the couple $(Y,X)$.
\end{prop}
\begin{Proof}
$X$ cannot be $\el^\infty$, neither $c_0$. So $X=\dis\bigcap_{p>\ga}\el^p$ for some $0\le\ga<+\infty$
or $X=\el^b$ for some $b\in(0,+\infty)$.

Let $X=\el^p$ with $0<b<+\infty$. We consider the space $Y_0=\dis\bigcap_{p>b}\el^p$. Obviously $Y\supset Y_0$.

According to Lemma \ref{lem2} there exists a vector space $F$ dense in $Y_0$ such that $F-\{0\}\subset Y_0\sm X\subset Y-Y_0$.
It remains yo show that $F$ is dense in $Y$.

Since convergence in $Y_0$ implies convergence in $Y$ and $F$ is dense in $Y_0$, it follows that the closure of $F$ in $Y$
contains $Y_0$. But $Y_0$ is dense in $Y$, because $c_{00}\subset Y_0$ is dense in $Y$. It follows that $F$ is dense in $Y$.
Thus, we have algebraic genericity in the case $X=\el^b$, $0<b<+\infty$. It remains the case
$X=\dis\bigcap_{p>\ga}\el^p$, $Y\neq X$ $X\subset Y\subset c_0$.

Then there exists $b\in(\ga,+\infty)$ such that $\el^b\subset Y$, $\el^b\neq Y$.

By the previous case, there exists a vector subspace $F$ of $Y$, dense in $Y$, such that $F-\{0\}\subset Y-\el^b$.

Since $\ga<b<+\infty$ it follows that $X=\dis\bigcap_{p>\ga}\el^p\subset\el^b$. This implies $Y-\el^b\subset Y-X$.
Thus, $F-\{0\}\subset Y-X$ and we have algebraic genericity in this case, as well.

The proof is completed. \qb
\end{Proof}
\section{Spaceability for the $\bbb{\el^p}$ spaces}\label{sec3}
\noindent

Let $X$ and $Y$ be two spaces as above that is, among $\el^\infty$, $c_0$, $\el^p$ and instersections
 of them $X\subset Y$, $X\neq Y$. We say that there is spaceability for the couple $(Y,X)$ if there exists a
 closed infinite dimensional subspace $F$ of $Y$ such that $F-\{0\}\subset Y\sm X$. We will show in this section
 that we always have spaceability.
\begin{prop}\label{prop4}
Let $Y=\el^\infty$ and $X\neq\el^\infty$ as above. Then we have spaceability for the couple $(\el^\infty,X)$.
\end{prop}
\begin{Proof}
Obviously $X\subset c_0$. Thus, it suffices to find a closed infinite dimensional vector subspace $F$
of $Y=\el^\infty$ so that $F-\{0\}\subset\el^\infty-c_0$. Let $A_m\subset N$, $m=1,2\ld$ be a sequence of
infinite subsets of $\{1,2,\ld\}$ pairwise disjoint. Let $F=\{f\in\el^\infty:f(A_m)$ is a singleton for each $m=1,2,\ld\}$.
It is easy to verify that $F$ is a closed infinite dimensional vector subspace of $\el^\infty$ and that $F-\{0\}\subset\el^\infty-c_0$.
The proof is\linebreak completed. \qb
\end{Proof}
\begin{prop}\label{prop5}
Let $X=\el^b$, $0<b<+\infty$ and $Y\supset X$, $Y\neq X$ as above. Then we have spaceability for the couple $(Y,\el^b)$.
\end{prop}
\begin{Proof}
Convergence $f^m\xrightarrow[m\ra+\infty]{}f$in $Y$ implies pointwise convergence $f^m(n)\underset{m}{\rightarrow}f(n)$
for all $n=1,2,\ld$.

Consider a sequence $A_j$, $j=1,2,\ld$ of infinite subset of $\{1,2,\ld\}$ pairwise disjoint. Let $y_j:\{1,2,\ld\}\ra\C$ be
such that $y_j(n)=0$ for all $n\in\{1,2,\ld\}\sm A_j$, $y_j\in Y$ but $y_j\notin\el^b$; that
is $\ssum_{n\in A_j}|y_j(n)|^b=+\infty$. Let $F$ be the closure in $Y$ of the linear span $\langle y_1,y_2,\ld\rangle$.
It is immediate that $F$ is a closed infinite dimensional vector subspace of $Y$. It remains to show that if $f\in F$, $f\neq0$,
than $f\notin\el^b$.

There exists a sequence $f^m\in\langle y_1,y_2,\ld\rangle$ with $f^m\xrightarrow[m\ra+\infty]{}f$ in $Y$. Since $f\neq0$,
there exists $n_0\in\{1,2,\ld\}$, such that $f(n_0)\neq0$. There exists $j_0\in\{1,2,\ld\}$ so that $n_0\in A_{j_0}$.
Otherwise $f^m(n_0)=0$ for all $m$, which implies that $f(n_0)=\dis\lim_m f^m(n_0)=0$. For each $m=1,2,\ld$ there
exists a constant $c_m$ such that $f^m|A_{j_0}=c_m\cdot y_{j_0}|_{A_{j_0}}$.

If $y_{j_0}(n_0)=0$, then $f^m(n_0)=0$ for all $m$. This implies that $f(n_0)=0$ which contradicts the fact that $f(n_0)\neq0$.
Therefore, $y_{j_0}(n_0)\neq0$.

Since $f(n_0)=\dis\lim_mf^m(n_0)=\dis\lim_{m}[c_my_{j_0}(n_0)]$, it follows that $\dis\lim_mc_m=\dfrac{f(n_0)}{y_{j_0}(n_0)}$
exists and is a constant $c\neq0$. It follows that for any $n\in A_{j_0}$ we
have $f(n)=\dis\lim_mf^m(n)=\dis\lim_mc_my_{j_0}(n)=cy_{j_0}(n)$. Therefore,
\[
\sum_{n\in A_{j_0}}|f(n)|^b=\sum_{n\in A_{j_0}}|c|^b|y_{j_0}(n)|^b=+\infty, \ \ \text{since} \ \ c\neq0.
\]
Thus, $f\notin\el^b$ and the proof is completed. \qb
\end{Proof}
\begin{prop}\label{prop6}
Let $X=\dis\bigcap_{p>a}\el^p$, $0\le a<+\infty$ and $Y\supset X$, $Y\neq X$ as above. Then we have spaceability for the
couple $(Y,X)$.
\end{prop}
\begin{Proof}
There exists $b\in(a,+\infty)$ such that $\el^b\subset Y$, $\el^b\neq Y$. Obviously, since $a<b$ we
have $X=\dis\bigcap_{p>a}\el^p\subset\el^b$. According to Proposition \ref{prop5}, there exists a closed
infinite dimensional vector subspace $F$ of $Y$, such that $F-\{0\}\subset Y-\el^b\subset Y-X$. Thus, we have
spaceability for the couple $(Y,X)$.

The proof is completed. \qb
\end{Proof}

Combining Propositions \ref{prop4}, \ref{prop5} and \ref{prop6} we obtain
\begin{thm}\label{thm7}
Let $Y,X$ as above, $Y\supset X$, $Y\neq X$, Then we have spaceability of the couple $(Y,X)$.
\end{thm}
\section{Continuation of the project}\label{sec4}
\noindent

Other chains of spaces, where we can investigate the same questions, are the followings:
\begin{enumerate}
\item[1)] $A(D)\subset H^\infty(D)\subset BM0A(D)\subset H^p(D)\subset H(D)$ and intersections of those spaces.
\item[2)] $A(D)\subset H^\infty(D)\subset$\, Bergman space\, $OL^p(D)\subset H(D)$.
\item[3)] For  $0\le p\le1$, $1\le\ga<+\infty$,\\
    $\el^p\subset\el^1=\Big\{\ssum^\infty_{n=0}a_nz^n,\ssum|a_n|<+\infty\Big\}\subset A(D)\subset H^\infty(D) \subset
    BMOA(D)\subset\\ \subset H^\ga(D)\subset c_0\subset\el^\infty\subset H(D)$.

This chain can also be continued as follows: $H^\ga(D)\subset H^\bi(D)\subset H(D)$ with $0<\bi<1$.
\item[4)] $0\le p\le2$, $2\le\de\le1$, $\el^p\subset\el^2=H^2(D)\subset H^\de(D)\subset H^1(D)\subset c_0\subset\el^\infty\subset H(D)$.
\item[5)] $0\le p\le2$, $0<\e<4$, $\el^p\subset\el^2=H^2(D)\subset OL^4(D)\subset OL^\e(D)\subset H(D)$.
\item[5$'$)] $H^p(D)\subset OL^{2p}(D)\subset H(D)$.
\item[6)] All previous with localized version of these spaces. As for example $H^1(D)\subset H^1_{(\al,\bi)}(D)$ localized version
\end{enumerate}
\section{Topological genericity for Hardy spaces}\label{sec5}
\noindent

In this section, as well as, in Section \ref{sec6}, we give two results in the frame of the project of Section \ref{sec4}
 and especially for Hardy spaces on the unit disc. \vspace*{0.5cm}\\
\noindent
{\bf Definition}. {\em Let $a<b$ and $p\in (0,+\infty)$. A holomorphic function $f$ on the open unit disc $D$ belongs
to the localised Hardy space $H^p_{[a,b]}(D)$, if $\dis\sup_{0<r<1}\dis\int^b_a|f(re^{i\vthi})|^pd\vthi<+\infty$.

If $b-a\ge2\pi$ then $H^p_{[a,b]}(D)$ coincides with the usual Hardy space $H^p=H^p(D)$. The spaces $H^p_{(a,b)}(D)$
are $F$-spaces under their natural topologies. A sequence $f_n\in H^p_{[a,b]}(D)$ converges in $H^p_{[a,b]}(D)$ to a
function $f\in H^p_{[a,b]}(D)$ iff $f_n\ra f$ uniformly on each compact subset of
 $D$ and $\dis\sup_{0<r<1}\dis\int^b_a|f_n(re^{i\vthi})-f(re^{i\vthi})|^pd\vthi\ra0$,
as $n\ra+\infty$}.

The space $H^p_{[a,b]}(D)$ contains $H^p(D)=H^p$ and the inclusion map is continuous.
\begin{thm}\label{thm8}
Let $p\in(0,+\infty)$. Then, the set $\{f\in H^p:$ for every $p<q<+\infty$ and every $a<b$ it holds
 $f\notin H^q_{[a,b]}(D)\}=\{f\in H^p:$ for all $p<q<+\infty$ and all
  $a<b$ it holds $\dis\sup_{0<r<1}\dis\int^b_a|f(re^{i\vthi})|^qd\vthi=+\infty\}$ is a $G_\de$ and dense
  subset of $H^p$. In particular it is non-void.
\end{thm}
\begin{Proof}
Fix $p<q<+\infty$ and $a<b$. We still show that the set
$A=\Big\{f\in H^p:\dis\sup_{0<r<1}\dis\int^b_a|f(re^{i\vthi})|^qd\vthi=+\infty\Big\}$ is a $G_\de$-dense
subset of $H^p$. Then varying $q$ to a sequence $q_n$ strictly decreasing and converging to $p$ and $a,b$ in
the set $Q$ of rational numbers, by denumerable intersection the result follows using Baire's theorem.

We have to show that the set $H^p\sm A$ can be written as a denumerable union of closed in $H^p$ sets with empty interiors.

We have $H^p\sm A=\dis\bigcup^\infty_{M=1}E_M$ with $E_M=\dis\bigcap_{0<r<1}E_{M,r}$ where
\[
E_{M,r}=\bigg\{f\in H^p:\dis\int^b_a|f(re^{i\vthi}|^qd\vthi\le M\bigg\} \ \ \text{and}
\]
\[
E_M=\bigg\{f\in H^p:\sup_{0<r<1}\int^b_a|f(re^{i\vthi})|^qd\vthi\le M\bigg\}.
\]
First we show that $E_{M,r}$ is closed in $H^p$, which implies that $E_M$ is also closed in $H^p$. Indeed,
let $f_n\in E_{M,r}$ converges to $f$ in the topology of $H^p$. This implies that $f_n$ converges to $f$
 uniformly on each compact subset of $D$. Thus, $\dis\int^b_a|f_n(re^{i\vthi})|^qd\vthi\ra\dis\int^b_a|f(re^{i\vthi})|^qd\thi$,
as $n\ra+\infty$. As $\dis\int^b_a|f_n(re^{i\vthi})|^qd\vthi\le M$ for all $n$, it follows that $f\in E_{M,r}$.
We have proved that the sets $E_{M,r}$ and $E_M$ are closed in $H^p$.

It remains to show that $E^\circ_M=\emptyset$. If not, pick $f\in E^\circ_M$. Since $E^\circ_M\subset E_M$, by
the definition of $E_M$, it follows that $f\in H^q_{[a,b]}(D)$.

Let $\oo=\dfrac{\al+\bi}{2}$. For an appropriate choice of $\ga>0$ the function $g(z)=\dfrac{1}{(z-e^{i\oo})^\ga}$
belongs to $H^p\sm H^q_{(a,b)}(D)$.

Since $H^p$ is a topological vector space the sequence $f+\dfrac{1}{n}g$ converges in the topology of $H^p$
towards $f$, as $n\ra+\infty$. Since $f$ belongs to the open in $H^p$ set $E^\circ_M$, it follows that there exists $n$ so that $f+\dfrac{1}{n}g\in E^\circ_M$. Since $f\in E^\circ_M$ and $E^\circ_M\subset E_M\subset H^q_{[n,b]}(D)$  it follows that $f,f+\dfrac{1}{n}g\in H^q_{[a,b]}(D)$. Since $H^q_{[a,b]}(D)$ is a vector space we conclude that $g$ belongs to $H^q_{[a,b]}(D)$. This is a contradiction and the proof is complete.  \qb
\end{Proof}
\section{Algebraic genericity in $\bbb{H^p}$ spaces}\label{sec6}
\noindent

In this section we prove a second original result in the frame of the project in \S\, 4.
\begin{thm}\label{thm9}
Let $0<p<q<+\infty$. Then $(H^p\sm H^q)\cup\{0\}$ contains a vector space dense in $H^p$.
\end{thm}
\begin{Proof}
There is $\ga>0$ so that the function $\dfrac{1}{(z-1)^\ga}$ belongs to $H^p-H^q$. Let $\oo_n=\dfrac{1}{n}$.
We consider the functions %
\[
f_n(z)=\frac{c_n}{(z-e^{i\oo_n})^\ga}, \ \ n=1,2,\ld\;.
\]
Choosing $c_n$ close to zero we obtain $d_p(f_n,0)<\dfrac{1}{n}$, where $d_p$ is the metric in $H^p$. Let $P_n$, $n=1,2,\ld$
be an enumeration of all polynomials with coefficients in $Q+iQ$. The sequence $P_n$, $n=1,2,\ld$ is dense in $H^p$.
Since $H^p$ does not contain isolated points, it follows that the sequence $P_n+f_n$, $n=1,2,\ld$ is also dense in $H^p$.
Thus, the linear space $F=\Big\langle f_n+P_n\Big\rangle^\infty_{n=1}$ is a vector space dense in $H^p$.

Let
\[
L=\la_1(f_1+P_1)+\cdots+\la_N(f_N+P_N), \ \  \la_N\neq0.
\]
We have to show that $L$ does not belong to $H^q$. Let $a<\oo_N<b$ with $\oo_1,\ld,\oo_{N-1}\notin[a,b]$. Then
\[
\sup_{0<r<1}\int^b_a|\la_k(f_k+p+k)(re^{i\vthi})|^qd\vthi<+\infty \ \ \text{for} \ \ k=1,\ld,N-1.
\]
and
\[
\sup_{0<r<1}\int^b_a|\la_N(f_N+P_N)(re^{i\vthi})|^qd\thi=+\infty.
\]
It follows that
\[
\sup_{0<r<1}\int^{2\pi}_0|L(re^{i\thi})|^qd\vthi=+\infty.
\]
Thus, $L$ does not belong to $H^q$ and the proof is complete. \qb
\end{Proof}

{\bf Acknowledgment:} We wish to thank R. Aron, L. Bernal\,-\,Gonzalez, I. Deliyanni and A. G. Siskakis for helpful
communications.

\vspace*{0.5cm}
\noindent
V. Nestoridis\\
National and Kapodistrian University of Athens\\
Department of Mathematics\\
Panepistemiopolis, 157 84\\
Athens,\\
Greece\\
e-mail: vnestor@math.uoa.gr \\


\begin{thebibliography}{99}
\bibitem{1} L. Bernal\,-\,Gonzalez, D. Pellegrino and J. B. Seoane\,-\,Sepulveda, Linear subsets of non linear sets
in topological vector spaces, Bulletin (New Series) of the AMS, vol. 51, No 1, January 2014, pages 71-130.
%
\bibitem{2} R. M. Aron, L. Bernal\,-\,Gonzalez, D. M. Pellegrino and J. B. Seone\,-\,Sepulveda, Lineability,
The Search for Linearity in Mathematics, Monographs and Research Notes in Mathematics, CRC Press, Taylor and Francis Group,
A chapman and Hall book.
    %
\bibitem{3} F. Bayart, K.-G. Grosse-Erdmann, V. Nestoridis and C. Papadimitropoulos, Abstract theory of universal series
and applications, Proc. London Math. Soc. 96 (2009) no 2, 417-468.
%
\bibitem{4} L. Bernal\,-\,Gonz\'{a}lez, H. J. Cabana\,-\,M\'{e}ndez, G. A. Mu\~{n}oz\,-\,Fern\'{a}der and Seoane\,-\,Sep\'{u}lveda,
Ordering among the topologies induced by various polynomial norms, Bull. Belgian Mathematical Society,
Simon Sterin 26(2019) no. 4, 481-492.
%
\bibitem{5} V. Nestoridis, A. G. Sirkakis, A. Stavrianidi, S. Vlachos. Generic non-extendability and total
unboundedness in function spaces, JMAA 475 (2019), no. 2 1720-1731, see also arxiv:1811.04408.
%
\bibitem{6} K. Kioulafa. On Hardy type spaces in strictly pseudoconvex domains and the density in these spaces,
 of certain classes of singular functions, JMAA (2020) no. 1, 123697, 18 pages, see also arxiv.1905.0476.
\end{thebibliography}
\end{document}